\documentclass[12pt]{article}
\usepackage{amsmath,amssymb,latexsym,mathrsfs,bm}
\usepackage{mathabx}
\usepackage[utf8]{inputenc}
\usepackage{tikz-cd}
\usetikzlibrary{cd}
\usepackage{comment}

\begin{document}
\title{Involutive Yang-Baxter groups never act as Frobenius groups}
\author{Arpan Kanrar$^{a\ast}$ and Wolfgang Rump$^{b}$  \\[2mm]                                  
{\small\it ${}^a$Harish-Chandra Research Institute}\\[-1.5mm]                                   
{\small\it A CI of Homi Bhabha National Institute}\\[-1.5mm]
{\small\it Chhatnag Rd, Jhunsi, Prayagraj - 211019, India}\\[1mm]                                
${}^b${\small\it Institute for
Algebra and Number Theory, University of Stuttgart}\\[-1.5mm]                                    
{\small\it Pfaffenwaldring 57, D-70550 Stuttgart, Germany}\\[3mm]
{\small\rm Dedicated to B.~V.~M.}\\[1mm]}

\date{}
\maketitle
\setcounter{footnote}{-1}
\renewcommand{\thefootnote}{\Roman{footnote}}
\footnotetext{\hspace*{-7mm} e-mail: {\tt arpankanrar@hri.res.in,
rump@mathematik.uni-stuttgart.de}\\
${}^\ast$Arpan Kanrar is supported by SRF-PhD fellowship from HRI.

\it \hspace*{-5mm} Key words and phrases. \rm Yang-Baxter Equation, cycle
set, braces, Frobenius action}
\thispagestyle{empty}

\newtheorem{prop}{Proposition}
\newtheorem{thm}{Theorem}
\newtheorem{lem}{Lemma}
\newtheorem{Definition}{Definition}
\renewcommand{\labelenumi}{\rm (\alph{enumi})}
\newcommand{\hs}{\hspace{2mm}}
\newcommand{\vsp}{\vspace{4ex}}
\newcommand{\vspc}{\vspace{-1ex}}
\newcommand{\hra}{\hookrightarrow}
\newcommand{\tra}{\twoheadrightarrow}
\newcommand{\md}{\mbox{-\bf mod}}
\newcommand{\Mod}{\mbox{-\bf Mod}}
\newcommand{\Mdd}{\mbox{\bf Mod}}
\newcommand{\mdd}{\mbox{\bf mod}}
\newcommand{\latt}{\mbox{-\bf lat}}
\newcommand{\Proj}{\mbox{\bf Proj}}
\newcommand{\Inj}{\mbox{\bf Inj}}
\newcommand{\Ab}{\mbox{\bf Ab}}
\newcommand{\CM}{\mbox{-\bf CM}}
\newcommand{\Prj}{\mbox{-\bf Proj}}
\newcommand{\prj}{\mbox{-\bf proj}}
\newcommand{\ra}{\rightarrow}
\newcommand{\eps}{\varepsilon}
\renewcommand{\epsilon}{\varepsilon}
\renewcommand{\phi}{\varphi}
\renewcommand{\rho}{\varrho}
\renewcommand{\hom}{\mbox{Hom}}
\newcommand{\ex}{\mbox{Ext}}
\newcommand{\rad}{\mbox{Rad}}
\renewcommand{\Im}{\mbox{Im}}
\newcommand{\oti}{\otimes}
\newcommand{\sig}{\sigma}
\newcommand{\en}{\mbox{End}}
\newcommand{\Qq}{\mbox{Q}}
\newcommand{\lra}{\longrightarrow}
\newcommand{\lras}{\mbox{ $\longrightarrow\raisebox{1mm}{\hspace{-6.5mm}$\sim$}\hspace{3mm}$}}
\newcommand{\Eq}{\Leftrightarrow}
\newcommand{\Equ}{\Longleftrightarrow}
\newcommand{\Ra}{\Rightarrow}
\newcommand{\Lra}{\Longrightarrow}
\newcommand{\A}{\mathbb{A}}
\newcommand{\N}{\mathbb{N}}
\newcommand{\Z}{\mathbb{Z}}
\newcommand{\C}{\mathbb{C}}
\newcommand{\Q}{\mathbb{Q}}
\newcommand{\R}{\mathbb{R}}
\renewcommand{\H}{\mathbb{H}}
\newcommand{\K}{\mathbb{K}}
\newcommand{\F}{\mathbb{F}}
\renewcommand{\P}{\mathbb{P}}
\newcommand{\B}{\mathbb{B}}
\newcommand{\rk}{\mbox{r}}
\newcommand{\p}{\mathfrak{p}}
\newcommand{\q}{\mathfrak{q}}
\renewcommand{\k}{\mathfrak{k}}
\newcommand{\AAA}{\mathfrak{A}}
\newcommand{\PPP}{\mathfrak{P}}
\renewcommand{\AA}{\mathscr{A}}
\renewcommand{\SS}{\mathscr{S}}
\newcommand{\DD}{\mathscr{D}}
\newcommand{\BB}{\mathscr{B}}
\newcommand{\LL}{\mathscr{L}}
\newcommand{\HH}{\mathscr{H}}
\newcommand{\NN}{\mathscr{N}}
\newcommand{\CC}{\mathscr{C}}
\newcommand{\EE}{\mathscr{E}}
\newcommand{\MM}{\mathscr{M}}
\newcommand{\OO}{\mathscr{O}}
\newcommand{\PP}{\mathscr{P}}
\newcommand{\II}{\mathscr{I}}
\newcommand{\TT}{\mathscr{T}}
\newcommand{\XX}{\mathscr{X}}
\newcommand{\YY}{\mathscr{Y}}
\newcommand{\FF}{\mathscr{F}}
\newcommand{\GG}{\mathscr{G}}
\newcommand{\RR}{\mathscr{R}}
\newcommand{\setm}{\smallsetminus}
\renewcommand{\le}{\leqslant}
\renewcommand{\ge}{\geqslant}
\newcommand{\rat}{\rightarrowtail}
\newcommand{\op}{^{\mbox{\scriptsize op}}}
\newcommand{\pf}{\it Proof.\hs\rm}
\newcommand{\bx}{\hspace*{\fill} $\Box$}
\newcommand{\ltra}{\lra\!\!\!\!\!\!\:\ra}
\newcommand{\subs}{\subset}
\newcommand{\sups}{\supset}
\newcommand{\subsn}{\subsetneq}
\newcommand{\nsubs}{\not\subset}
\newcommand{\pt}{\makebox[0pt][r]{\bf .\hspace{.5mm}}}
\newcommand{\dpt}{\makebox[0mm][l]{.}}
\newcommand{\dpc}{\makebox[0mm][l]{,}}
\newcommand{\noth}{\varnothing}
\newcommand{\da}{\downarrow\!\!}
\newcommand{\ua}{\uparrow\!\!}

{\small\bf Abstract.} {\small A conjecture of S. Ram\'{\i}rez states that every
indecomposable non-degenerate involutive set-theoretic solution to the Yang-Baxter
equation with dihedral permutation group of order $2n$ has cardinality $2n$. The
conjecture is verified for odd $n$ and disproved for even $n$. The proof for odd
$n$ is obtained from the more general result that the permutation group of a finite
solution never acts as a Frobenius group.}

\section{Introduction}

After Drinfeld's suggestion \cite{Dri} to study solutions $r$ to the set-theoretic
Yang-Baxter equation
\begin{equation}\label{1}
(r\times 1_X)(1_X\times r)(r\times 1_X)=(1_X\times r)(r\times 1_X)(1_X\times r),  
\end{equation}
solutions were found on symplectic manifolds \cite{WX}, and then in connection with
Artin-Schelter regular rings \cite{GI8, GI-V} and algebras of I-type \cite{TV}. 
A solution $r(x,y)=({}^xy,x^y)$ of Eq.~(\ref{1}) is said to be {\em non-degenerate}
\cite{ESS} if its component maps $y\mapsto{}^xy$ and $x\mapsto x^y$ are bijective.
If $X$ does not admit a non-trivial partition $X=Y\sqcup Z$ with induced solutions on
$Y$ and $Z$, the solution $(X,r)$ is said to be {\em indecomposable}. A systematic
study of non-degenerate {\em involutive} solutions $r=r^2$ was initiated by Etingof
et al. \cite{ESS} who associate a {\em structure group} $G_X$ to any solution of
Eq.~(\ref{1}).
For a finite non-degenerate involutive solution on $X$, a finite factor group $G(X)$
of $G_X$, acting on $X$
by permutations, has become an important invariant of the solution. For example, 
such a solution $(X,r)$ is indecomposable if and only if the {\em permutation group}
$G(X)$ acts transitively on $X$. 

Recently, Ram\'{\i}rez \cite{Ram} conjectured that a dihedral permutation group
$G(X)$ of an indecomposable non-degenerate involutive solution of Eq.~(\ref{1}) must
be of
cardinality $|X|$. If true, Ram\'{\i}rez' conjecture would provide a highly
remarkable, exceptional family of involutive Yang-Baxter groups \cite{CJR} and
corresponding solutions to Eq.~(\ref{1}). For example, such solutions would be
their own universal covering \cite{Ind} (see Section~3 for details).   

In this paper, we verify the conjecture for dihedral groups $D_{2n}$ with odd $n$,
and provide a simple counter-example for even $n$. Dihedral groups
$D_{2n}$ with odd $n$ are Frobenius groups. Ram\'{\i}rez \cite{Ram} proved that
the conjecture is ``almost true'' in the sense that a possible counter-example
must satisfy $|X|=n$ whenever $G(X)=D_{2n}$. We show that in the latter case,
$D_{2n}$ acts as a Frobenius group on $X$. To exclude this possibility, we prove
that more generally, the permutation group of a finite non-degenerate involutive
solution $(X,r)$ of Eq.~(\ref{1}) never acts as a Frobenius group on $X$.

The occasion to write this paper arose from a speculation on a question of
Ram\'{\i}rez and Vendramin \cite{RaV} which asks whether a finite
non-degenerate involutive solution $(X,r)$ is decomposable if for some $x\in X$,
the permutation $\sigma(x)$ contains a cycle of length coprime to $|X|$.
Castelli \cite{Cas} verified this for {\em multipermutation} solutions (see
\cite{ESS} or Section~2), while the first
author obtained a positive answer if $G(X)$ is nilpotent \cite{Kan}. If the
decomposibility
question \cite{RaV} is true in general, this would give an alternative proof of
our positive result for dihedral groups $D_{2n}$ with $n$ odd.

\section{Cycle sets and braces} 

In this section we recall some facts about cycle sets and braces as far as needed
below. For a finite non-degenerate involutive solution $r(x,y)=({}^xy,x^y)$ and
$x\in X$, let $\sigma(x)\colon X\ra X$ be the self-map of $X$ given by
$\sigma(x)^{-1}(y)=y^x$. The binary operation $x\cdot y:=\sigma(x)(y)$ then satisfies
\begin{equation}\label{2}
(x\cdot y)\cdot(x\cdot z)=(y\cdot x)\cdot(y\cdot z)  
\end{equation}
for all $x,y,z\in X$. A set with a binary operation $(X;\cdot)$ satisfying
Eq.~(\ref{2}) where the maps $\sigma(x)\colon X\ra X$ with $\sigma(x)(y)=x\cdot y$
are bijective is said to be a {\em cycle set} \cite{YBE}. With $x^y:=\sigma(y)^{-1}(x)$
as above, every cycle set $X$ determines a solution to Eq.~(\ref{1}):
\begin{equation}\label{3}
r(x,y)=({}^xy,x^y)  
\end{equation} 
with
\begin{equation}\label{4}
{}^xy:=x^y\cdot y.  
\end{equation}
By \cite{YBE}, Proposition~1 and Theorem~2, this gives a one-to-one correspondence
between finite cycle sets and finite
non-degenerate involutive set-theoretic solutions to the Yang-Baxter equation.
Being given by a single operation, cycle sets are much easier to handle than solutions
$r\colon X\times X\ra X\times X$ to Eq.~(\ref{1}). The correspondence extends to
infinite solutions if a non-degeneracy condition for cycle sets is assumed, which states    
the square map $T\colon X\ra X$ with $T(x):=x\cdot x$ is bijective. For finite cycle
sets, it is important to know that this map is always bijective. The somewhat asymmetric
equations (\ref{3}) and (\ref{4}) can be put into a more appealing form:
$$r(y\cdot x,y)=(x\cdot y,x).$$

By \cite{YBE}, Proposition~5, the operation of a cycle set $X$ extends in a unique
fashion to the free abelian monoid $\N^{(X)}$, so that the following equations hold
for $a,b,c\in\N^{(X)}$:
\begin{align}
a\cdot(b+c) &= (a\cdot b)+(a\cdot c) \label{5}\\
(a+b)\cdot c &= (a\cdot b)\:\cdot\:(a\cdot c). \label{6}
\end{align}
If $X$ is non-degenerate (e.~g., finite), the unique extension further extends to
the free abelian group $\Z^{(X)}$ (\cite{YBE}, Proposition~6). An abelian group
$A$ with a binary operation $(A;\cdot)$ satisfying Eqs.~(\ref{5}) and (\ref{6}) is
said to be a {\em linear cycle set} \cite{YBE} or a {\em brace} \cite{Bra}. By
Eq.~(\ref{6}), every brace is a cycle set. Conversely, every
non-degenerate cycle set $X$ embeds into a brace $A_X$ with additive group
$\Z^{(X)}$.

There are more operations on a brace $A$. For example, there is a second
group operation
\begin{equation}\label{7}
a\circ b:=a^b+b  
\end{equation} 
from which the cycle set operation can be recovered. Guarnieri and Vendramin \cite{GV}
showed that a set $(A,+,\circ)$ with two group operations, $(A;+)$ being commutative,
is a brace if and only if the single equation
$$(a+b)\circ c=(a\circ c)-c+(b\circ c)$$
holds in $A$. (This led to the concept of {\em skew-brace} where the commutativity of
$(A;+)$ is dropped.) The group $(A;\circ)$ of a brace $A$ acts from the left via
$\sigma(a)$ on the additive group by Eq.~(\ref{5}) and
\begin{equation}\label{8}
(a\circ b)\cdot c=a\cdot(b\cdot c),  
\end{equation} 
so that $A$ becomes a left $A^\circ$-module. Similarly, $A$ is a right $A^\circ$-module
under the exponential action:
\begin{align*}
(a+b)^c &= a^c+b^c\\ 
a^{b\circ c} &= (a^b)^c. 
\end{align*} 
The operation (\ref{7}) generalizes Jacobson's circle operation
\begin{equation}\label{9}
a\circ b:=ab+a+b  
\end{equation} 
which makes the Jacobson radical of any (associative unital) ring
into a brace. Following Jacobson's terminology for radical rings \cite{Jac}, the
group $A^\circ:=(A;\circ)$ of a brace $A$ is called the {\em adjoint group}, and
the inverse of an element $a\in A$ is denoted by $a'$. Note that the neutral element
of the adjoint group is 0, as in the adjoint group of a radical ring. 

The ring operation $ab$ (juxtaposition) in Eq.~(\ref{9}) is useful also for
braces. For example, a {\em right ideal} \cite{Bra} of a brace $A$ is defined to be
an additive subgroup $I$ which is closed with respect to right multiplication: $a\in
I$ and $b\in A$ implies $ab\in I$. If also $ba\in I$, then $I$ is called an {\em ideal}
\cite{Bra}. In terms of the above operations, an additive subgroup $I$ of $A$ is a
{\em right ideal} if $I$ is an $A^\circ$-submodule, and an {\em ideal} if, in addition,
$I$ is a
normal subgroup of $A^\circ$. In contrast to ring theory, left ideals play no significant
part in the theory of braces. For an ideal $I$ of a brace $A$, the residue classes
$a+I$ form a brace $A/I$, and there is a canonical brace morphism $A\tra A/I$, like in
ring theory. 

An important ideal of any brace $A$ is the {\em socle}
\begin{equation}\label{10}
\mbox{Soc}(A):=\{ a\in A\:|\:\forall\,b\in A\colon a\cdot b=b\}.  
\end{equation} 
For a non-degenerate cycle set $X$, the {\em permutation group} $G(X)$, generated by the
$\sigma(x)$ with $x\in X$, is the adjoint group
of $A(X):=A_X/\mbox{Soc}(A_X)$. Indeed, Eq.~(\ref{10}) shows that the kernel of the brace
morphism $A_X\tra A_X/\mbox{Soc}(A_X)$ consists of the elements for which the adjoint
action on $X$ is trivial. 

In general, the cycle set morphism
\begin{equation}\label{11}                                                                       
\sigma\colon X\ra A(X)                                                                           
\end{equation}
is not injective. Its image $\sigma X$ is called the {\em retraction} of $X$. A
non-degenerate cycle set $X$ is said to be {\em irretractable} if $\sigma$ is injective,
otherwise {\em retractable}. If
some iterate $\sigma^n X$ is a singleton, the solution to Eq.~(\ref{1}) associated with
$X$ is said to be a {\em multipermutation} solution \cite{ESS}.

The brace $A_X$ is functorial in $X$:
\begin{prop}\pt\label{p1}
Let $f\colon X\ra Y$ be a morphism of non-degenerate cycle sets. There is
a unique brace morphism $A_f\colon A_X\ra A_Y$ such that the following diagram
commutes:
\begin{center}
	\begin{tikzcd}
		X\arrow[d,hook]\arrow[r,"f"]& Y\arrow[d,hook]\\
		A_X\arrow[r,"A_f"] &A_Y.
	\end{tikzcd}
\end{center}
\end{prop}

\pf As the additive group of $A_X$ is free abelian, $f$ extends uniquely to an additive
group homomorphism $A_f$ so that the diagram commutes. For simplicity, let us
write $f$ instead of $A_f$ in the following argument. By Eq.~(\ref{8}) and induction,
the sub-cycle set $X$ of $A_X$ is invariant under the adjoint action of $G_X$. For
$a\in A_X$ and $x,y\in X$, it follows that Eq.~(\ref{6}) gives  
\begin{gather*}
f\bigl((a+x)\cdot y\bigr)=f\bigl((a\cdot x)\cdot(a\cdot y)\bigr)=f(a\cdot x)\cdot
f(a\cdot y)\\  
f(a+x)\cdot f(y)=\bigl(f(a)+f(x)\bigr)\cdot f(y)=\bigl((f(a)\cdot f(x)\bigr)\cdot
\bigl((f(a)\cdot f(y)\bigr).
\end{gather*} 
Thus, if $f(a\cdot z)=f(a)\cdot f(z)$ holds for all $z\in X$, then $f\bigl((a+x)
\cdot y\bigr)=f(a+x)\cdot f(y)$. Conversely, assume that $f\bigl((a+x)                  
\cdot y\bigr)=f(a+x)\cdot f(y)$ holds for all $y\in X$, this implies that
$f(a\cdot x)\cdot f(a\cdot x)=\bigl((f(a)\cdot f(x)\bigr)\cdot\bigl((f(a)\cdot f(x)
\bigr)$. Since $X$ is non-degenerate, it follows that $f(a\cdot x)=f(a)\cdot f(x)$. 

\vspc
So the equation $f(a\cdot x)=f(a)\cdot f(x)$ (with $a\in A_X$ and $x\in X$) remains 
valid if an element of $X$ is added to or subtracted from $a$. As the equation holds 
for $a\in X$, it does hold in general. Now Eq.~(\ref{5}) shows that $x$ in the equation 
can be replaced by any element of $A_X$. Thus $A_f$ is a morphism of linear
cycle sets, hence a morphism of braces. \bx

Factoring out the socle, one should expect that $A(X)$ is functorial as well.
However, the following example shows that this is not the case if $f$ is not
surjective:

\noindent {\bf Example 1.} Let $X$ be an irretractable cycle set with an element
$x=x\cdot x$. Assume that $\sigma(y)\ne 1_X$ for all $y\in X$. A 4-element
cycle set $X$ of that type is given in Example~2 of the next section. Regard $X$ as
a sub-cycle set of $A_X$. The sub-cycle set $Y:=X\sqcup\{ 0\}$ of $A_X$ is still
irretractable, and $Z:=\{x,0\}$ is a retractable sub-cycle set of $Y$. As the composed
cycle set morphism $Z\hra Y\stackrel{\sigma}{\lra} A(Y)$ is injective, it cannot
factor through $\sigma\colon Z\ra A(Z)$. Thus $A_X$ cannot be replaced by $A(X)$ in
the commutative diagram of Proposition~\ref{p1}.

For surjective cycle set morphisms, we have the following corollary:

\vspace{3mm}
\noindent {\bf Corollary.} \it Let $f\colon X\tra Y$ be a surjective morphism of
non-degenerate cycle sets. There is a unique brace morphism $A(f)\colon A(X)\tra A(Y)$
such that the following diagram commutes:
\begin{center}
	\begin{tikzcd}
		X\arrow[d,"\sigma"]\arrow[r,"f",twoheadrightarrow]& Y\arrow[d,"\sigma"]\\
		A(X)\arrow[r,"A(f)",twoheadrightarrow] &A(Y).
	\end{tikzcd}
\end{center}

\vspace{3mm}
\pf We only have to show that $A_f\colon A_X\ra A_Y$ maps $\mbox{Soc}(A_X)$ into
$\mbox{Soc}(A_Y)$. Thus, let $a\in\mbox{Soc}(A_X)$ be given. For all $x\in X$, we have
$A_f(a)\cdot A_f(x)=A_f(a\cdot x)=A_f(x)$. Since $f$ is surjective and $Y=A_f(X)$
generates $A_Y$, this implies that $A_f(a)\in\mbox{Soc}(A_Y)$. Whence $A_f$
induces a brace morphism $A(f)\colon A(X)\ra A(Y)$. Since $\sigma Y$ generates
$A(Y)$, the brace morphism $A(f)$ is surjective. \bx

For a finite brace $A$, let $\pi(A)$ be the set of primes dividing the order of $A$. The
primary components $A_p$ of the additive group are right ideals, and we have a
decomposition
\begin{equation}\label{12}
A=\bigoplus_{p\in\pi(A)}A_p.   
\end{equation} 

The following result, which sheds some light upon the Yang-Baxter equation, is part of
\cite{Cyc2}, Proposition~2. For completeness, we include a proof.
\begin{prop}\pt\label{p2}
Let $A=I\oplus J$ be a brace with right ideals $I$ and $J$. Up to isomorphism, the
brace $A$ is uniquely determined by the sub-braces $I$ and $J$. For $a\in I$ and $b\in
J$, the mutual adjoint actions between $I$ and $J$ are determined by
\begin{equation}\label{13}
b\circ a={}^ba\circ b^a,  
\end{equation} 
where as in Eq.~$(\ref{4})$, ${}^ba:=b^a\cdot a$.
\end{prop}

\pf If $I$ and $J$ are given, the brace structure of $A$ is determined by the mutual
adjoint actions between $I$ and $J$. Thus, it is enough to verify Eq.~(\ref{13}). Note
that ${}^ba\in I$ and $b^a\in J$. By Eq.~(\ref{7}), we have $b\circ a=b^a+a=a+b^a=
(b^a\cdot a)^{b^a}+b^a=(b^a\cdot a)\circ b^a={}^ba\circ b^a$. \bx

Recall that a sub-cycle set $X$ of a brace $A$ is said to be a
{\em cycle base} \cite{Bra} if $X$ is $A^\circ$-invariant and generates the additive
group of $A$. For example,
the retraction $\sigma X$ of a non-degenerate cycle set $X$ is a cycle base in $A(X)$.

\vspace{3mm}
\noindent {\bf Corollary.} \it Let $A=I\oplus J$ be a brace with an ideal $I$ and a
right ideal $J$. For $a\in I$ and $b\in J$, the equations $a\cdot b=b$ and $b\cdot a=
b\circ a\circ b'$ hold. \rm 

\vspace{3mm}
\pf Replacing $b$ by $a\cdot b$ in Eq.~(\ref{13}) yields $(a\cdot b)\circ a=(b\cdot a)
\circ b$. Thus $(a\cdot b)\circ b'=(b\cdot a)\circ b\circ a'\circ b'$. Since $I^\circ$
is normal in $A^\circ$, this implies that $(a\cdot b)\circ b'\in I\cap J=0$. Whence
$a\cdot b=b$, and thus $b\cdot a=b\circ a\circ b'$. \bx 

The following observation is trivial, but fundamental.  
\begin{prop}\pt\label{p3}
Let $X$ be a non-degenerate cycle set. The map $(\ref{11})$ is equivariant with 
respect to the adjoint group $G(X)$. 
\end{prop} 

\pf Every element of $G(X)$ is a finite product $g_1\circ\cdots\circ g_n$ with $g_i=
\sigma(x_i)$ or $g_i=\sigma(x_i)'$ for some $x_i\in X$. By induction, we can assume 
that $n=1$. Thus, in the first case, we have to verify that $\sigma
\bigl(\sigma(x)(y)\bigr)=\sigma(x)\cdot\sigma(y)$. This follows since (\ref{11}) is
a cycle set morphism. Hence $\sigma(y)=\sigma(x)'\cdot\sigma\bigl(\sigma(x)(y)\bigr)$.
Replacing $y$ by $\sigma(x)'(y)$, we get $\sigma\bigl(\sigma(x)'(y)\bigr)=\sigma(x)'
\cdot\sigma(y)$. \bx  

\section{Frobenius action and Ramirez' conjecture} 

Let $G$ be a transitive permutation group, acting on a finite set $X$. Assume that
the elements $g\ne 1$ in $G$ leave at most one $x\in X$ fixed. If some $g\ne 1$
has a fixed point in $X$, then $G$ is called a {\em Frobenius group}. Frobenius'
theorem \cite{Frob} states that the elements which leave no element
of $X$ fixed together with 1 form a normal subgroup $K$ of $G$, which is now called
the {\em Frobenius kernel}. Thus $G=K\rtimes G_x$ for any stabilizer $G_x$ of some
$x\in X$. Therefore $G_x$ is said to be a {\em Frobenius complement} of the action. 
To this date, a proof of Frobenius' theorem has not been found without character
theory. 

A simple calculation (\cite{Hup}, V.8.3) shows that $|H|$ divides $|K|-1$. In 
particular, the orders of $H$ and $K$ are relatively prime. Another remarkable fact
about Frobenius groups is that the Frobenius kernel and the set of Frobenius 
complements are unique (\cite{Hup}, V.8.17).     

In this section, we show for a finite cycle set $X$, that the action of the
permutation group $G(X)$ on $X$ cannot be a Frobenius action. Note that by
\cite{ESS}, Theorem~2.15, $G(X)$ is solvable. If the Frobenius complements are
solvable, Shaw \cite{Sha} found a character-free proof of Frobenius' theorem 
using the transfer map.

\begin{thm}\pt\label{t1}
Let $X$ be a finite cycle set. Then the permutation group $G(X)$ does not act as a
Frobenius group on $X$.  
\end{thm} 

\pf Suppose that $G(X)$ acts as a Frobenius group on $X$. Thus $G(X)$ is the adjoint
group of a brace $A(X)$. Let $K$ be the Frobenius kernel of the action. So the order
of $K$ is relatively prime to each Frobenius complement. Hence Eq.~(\ref{12}) gives
a decomposition $A(X)=I\oplus H$ into right ideals with $|I|=|K|$. By Hall's
theorem (\cite{Rob}, 9.1.7), $K$ is conjugate to $I$. Since
$K$ is a normal subgroup of $G(X)$, it follows that $K=I$, and thus $K$ is a brace
ideal of $A(X)$. Thus $A(X)=K\rtimes H$ is a semidirect product of braces, and
$H$ is a Frobenius complement. 

\vspc
Let $x\in X$ be the fixed point of $H$. By Proposition~\ref{p3}, the map (\ref{11})
is $G(X)$-equivariant. Thus $\sigma(x)=k+h$ with $k\in K$ and $h\in H$, and
$h\cdot(k+h)=h\cdot\sigma(x)=\sigma(h(x))=\sigma(x)=k+h$. In particular,
$h\cdot k=k$, that is, $h\circ k\circ h'=k$ by the corollary of Proposition~\ref{p2}.
So $h$ and $k$ commute. For the action of $k$ on $X$,
this implies that $k(x)=k(h(x))=h(k(x))$. So $h$ fixes $x$ and $k(x)$, which
gives $k(x)=x$. Whence $k=0$. 

\vspc
Now the Frobenius complement $H$ fixes $x$. Hence $H$ fixes $\sigma(x)=h$. On the other
hand, the corollary of Proposition~\ref{p2}, applied to the semidirect product $A(X)=
K\rtimes H$ of
braces, implies that $K^\circ$ acts trivially on $H$. So $G(X)$ fixes $\sigma(x)$.
Since $\sigma X$ is a cycle base of $A(X)$, it follows that $\sigma(x)$ generates
$A(X)$. Whence $A(X)\subs H$, a contradiction. \bx

Now let us turn to Ram\'{\i}rez' conjecture \cite{Ram}. By the corollary of
Proposition~\ref{p1}, a surjective morphism $f\colon Y\tra X$ of non-degenerate
indecomposable cycle sets induces a surjective brace morphism $A(f)\colon A(Y)\tra
A(X)$. If $A(f)$ is bijective, the morphism $f$ is said to be a {\em covering}
\cite{Ind} of $X$. Choose a fixed $x\in X$. By \cite{Ind}, Corollary~4.6, there is a
one-to-one correspondence, up to isomorphism, between the coverings of $X$ and the
subgroups of the {\em fundamental group} $\pi_1(X,x):=\{ a\in G(X)\:|\:a(x)=x\}$ of
$X$. Like in topology, there is a {\em universal covering} $p\colon\widetilde{X}\tra
X$ such that $\widetilde{X}$ has trivial fundamental group, that is $G(\widetilde{X})$
acts simply transitively on $\widetilde{X}$. Note that by the definition of a
covering, $G(\widetilde{X})$ is isomorphic to $G(X)$. So there is a bijection between
$\widetilde{X}$ and $G(X)$, which endows $G(X)$ with a second cycle set structure
coming from identification with $\widetilde{X}$. 

Ram\'{\i}rez' conjecture states that a finite indecomposable cycle set with dihedral
permutation group
$D_{2n}$ has trivial fundamental group, i.~e. it coincides with its universal covering.
The following example restricts the range of validity of the conjecture to odd $n$.

\noindent {\bf Example 2.} Let $X=\{ 1,2,3,4\}$ be the cycle set given by the permutations
$$\sigma(1)=(2,4),\quad\sigma(2)=(1,2,3,4),\quad\sigma(3)=(1,4,3,2),\quad\sigma(4)=(1,3).$$
This cycle set occurs as Example~14 in \cite{CCP}. Its adjoint group is dihedral of
order 8. The brace $A(X)$ is listed as $B_7$ in \cite{finbra}, and is given explicitly in
Example~2 of \cite{Bra}.

Thus, it remains to consider the odd case. It is well known that dihedral groups $D_{2n}$
with odd $n$ are Frobenius groups. By lack of a reference, we give a short proof.

\begin{prop}\pt\label{p4}
For odd $n\ge 3$, the dihedral group $D_{2n}$ is a Frobenius group.     
\end{prop}

\pf $D_{2n}$ is a Coxeter group, generated by two reflections $b,c$, with relations
$b^2=c^2=1$ and $(bc)^n=1$. Alternatively, $D_{2n}$ is generated by the rotation $a:=bc$
and the reflection $b$, so that the relations become $a^n=1$ and $bab=a^{-1}$. (Indeed,
$bab=b(bc)b=cb=(bc)^{-1}$.) There is a cyclic subgroup $\langle a\rangle$, and the
remaining $n$ elements are reflections. If $n$ is odd, each subgroup generated by a
reflection is a Sylow subgroup. All these 2-element subgroups are conjugate by Sylow's
theorem. Since $n$ is odd, each of them coincides with its normalizer. Thus $D_{2n}$ is
a Frobenius group with Frobenius kernel $\langle a\rangle$ and Frobenius complement
$\langle b\rangle$. \bx

For dihedral groups $D_{2n}$ with odd $n$, Ram\'{\i}rez' conjecture can now be verified.
\begin{thm}\pt\label{t2}
Let $D_{2n}$ be a dihedral group with odd $n\ge 3$. Every indecomposable cycle set $X$
with permutation group $D_{2n}$ is of cardinality $2n$.   
\end{thm}

\pf Since $X$ is indecomposable, $G(X)$ acts transitively on $X$. Let $H$ be the stabilizer
of some $x\in X$. So the other elements of $X$ are fixed by some conjugate of $H$. Since $G(X)$
acts faithfully on $X$, the intersection of all conjugates of $H$ is trivial. Thus $H$
does not contain a rotation of $D_{2n}$. As the product of two reflections is a rotation,
$H$ contains at most one reflection, which yields $|H|\le 2$. If $|H|=2$, the
action of $G(X)$ on $X$ is equivalent to the action by conjugation on its reflections.
Indeed, $gHg^{-1}\mapsto g(x)$ with $g\in G(X)$ is a $G(X)$-equivariant isomorphism
between the $G(X)$-set of reflections in $G(X)$ (under conjugation) and the $G(X)$-set $X$.
By Proposition~\ref{p4}, $G(X)$ is a Frobenius group, which is impossible by Theorem~\ref{t1}.  
Whence $|H|=1$ and $|X|=2n$. \bx

\noindent {\bf Remark.} Ram\'{\i}rez and Vendramin (\cite{RaV}, Question 3.16) asked
whether a finite cycle set $X$ must be decomposable if some $\sigma(x)$ with $x\in X$
contains a cycle
of length coprime to $|X|$. If true, we could get a simpler proof of Theorem~\ref{t1}:
Since $\sigma X$ is a cycle base of $A(X)$, it cannot be contained in the Frobenius
kernel. Hence it meets a Frobenius complement $H$, that is, a stabilizer of some $x\in X$.
As the order of $H$ is coprime to $|X|$, the cycle set $X$ would be decomposable, a
contradiction.

\noindent {\bf Acknowledgement.} The first author is thankful to Prof. Manoj Kumar
Yadav for a discussion on Frobenius groups.

\end{document}